\documentclass[letterpaper, 10 pt, conference]{ieeeconf}  

\IEEEoverridecommandlockouts                              
\usepackage{graphics} 
\usepackage{graphicx}
\usepackage{epsfig} 
\usepackage{subfigure}
\usepackage{amsmath} 
\usepackage{amssymb}  
\usepackage{latexsym}
\usepackage{url}
\usepackage{color}
\usepackage{xcolor}
\usepackage{subfigure}

\newtheorem{defi}{Definition}
\newtheorem{prop}{Proposition}
\newtheorem{theorem}{Theorem}
\newtheorem{assumption}{Assumption}
\newtheorem{lemma}{Lemma}
\newtheorem{remark}{Remark}
\newtheorem{corollary}{Corollary}

\title{\LARGE \bf
Impulsive Control for G-AIMD Dynamics with\\ Relaxed and Hard Constraints
}

\author{Konstantin Avrachenkov, Alexei Piunovskiy and Yi Zhang
\thanks{K. Avrachenkov is with Inria Sophia Antipolis, 2004 Route des Lucioles, 06902,
Sophia Antipolis, France, {\tt k.avrachenkov@inria.fr}}%
\thanks{A. Piunovskiy and Y. Zhang are with the Department of Mathematical Sciences,
University of Liverpool, UK, {\tt piunov@liverpool.ac.uk}, {\tt yi.zhang@liverpool.ac.uk}}%
}

\begin{document}

\maketitle
\thispagestyle{empty}
\pagestyle{empty}

\begin{abstract}
Motivated by various applications from Internet congestion control to power
control in smart grids and electric vehicle charging, we study Generalized Additive Increase Multiplicative Decrease
(G-AIMD) dynamics under impulsive control in continuous time with the time average alpha-fairness criterion.
We first show that the control under relaxed constraints can be described by a threshold.
Then, we propose a Whittle-type index heuristic for the hard constraint problem.
We prove that in the homogeneous case the index policy is asymptotically optimal
when the number of users is large.
\end{abstract}


\section{Introduction}

For nearly two decades Additive Increase Multiple Decrease (AIMD) mechanism was one of the main
components in the TCP/IP protocol regulating data traffic across the Internet \cite{S93}.
In the absence of significant queueing delay, AIMD increases the data sending rate
linearly in time until packet loss and then drastically, in a multiplicative fashion,
reduces the sending rate. However, in the most recent versions of TCP (Compound \cite{Tetal06} in Windows
and Cubic \cite{HRX08} in Linux), the linear growth function has been changed to non-linear
functions to enable agile adaptation of the data sending rate. Such modifications can
be viewed as particular cases of Non-linear AIMD (NAIMD) dynamics. The possibilities
of non-linear modifications of AIMD are really endless. A thorough classification of NAIMD dynamics,
together with the analysis of some NAIMD classes, can be found in the book \cite{AIMDbook}. Here we consider
one important class of NAIMD dynamics, which we refer to as Generalized AIMD (G-AIMD) \cite{Aetal15}.
In the G-AIMD dynamics the acceleration of the sending rate in the increase phase depends
on the current value of the rate.

The other important recent development in the Internet architecture is the introduction of Software-Defined
Networking (SDN) technology \cite{Ketal15}. The SDN technology allows much finer control of resource
allocation (e.g., bandwidth allocation) in a network. Motivated by this opportunity, in the
present work we study the control of G-AIMD dynamics. In the networking context, when
allocating resource, it is very common to use some fairness function as optimization
objective. In the foundational work \cite{KMT98}, the authors proposed to use proportional
fairness in the context of the network utility maximization problem. Then, in \cite{MW00}
the $\alpha$-fairness function was proposed, which generalizes the proportional fairness
and gives max-min fairness and delay fairness as the other important particular cases.
A very good review of the network utility maximization problem can be found in \cite{S12}.
Most of the works on the resource allocation problem concern with long-term fairness,
which ignores instantaneous oscillations of the sending rate or short-term fairness.
Short-term fairness is particularly important in wireless and electrical networks.
Following \cite{AAR12}, in this work we optimize the integral of the $\alpha$-fairness function
over time, which represents short-term fairness.

We would like to note that recently AIMD and more generally NAIMD found new applications in
smart electrical grids \cite{CLRS14,HW13,KR11} and in power control for charging
electric vehicles \cite{AIMDbook,ElecVehicBook,Netal17}. We hope that our findings will also
be useful in these application domains.

Let us specifically describe our contributions: in the next section we formulate
the problem of short-term $\alpha$-fairness for resource allocation among G-AIMD
users as an impulsive control problem under constraints with time average criterion. We would
like to note that our impulsive control is different from the standard impulsive
control setting \cite{MR12}, where there is a constraint on the number of impulses
or on the total variation of the impulse control. In our case, we have
only a constraint on the system state. The present work also represents an advance with respect
to our previous work \cite{Aetal15}, where we have not considered the setting with constraints.
Here we consider both hard and relaxed constraints. In Section~\ref{sec:relaxed} we show that
in the case of the relaxed constraints, the optimal impulsive control of the G-AIMD
dynamics can be given in the threshold form. Then, in Section~\ref{sec:hard} we propose
a heuristic, which is similar in spirit to the celebrated Whittle index \cite{W88}.
We would like to note that in the past several attempts to prove indexability
of AIMD \cite{AADJ13,JS12} and G-AIMD \cite{ABP17} dynamics have been made. However,
to the best of our knowledge, it is for the first time that we prove the indexability of the
G-AIMD dynamics without any artificial conditions. We were able to make this
theoretical advance largely thanks to the framework of impulsive control in
the continuous time. The previous works on TCP indexability are all in discrete
time, and some are also in the discrete state space but \cite{ABP17} is in the continuous
state space.
Similarly to \cite{WW90}, we are able to show that in the homogeneous case the index policy
is asymptotic optimal in the regime of a large number of users. As a by-product, we prove
the global stability of the AIMD dynamics and the local asymptotic stability
of the G-AIMD dynamics under the index policy in the homogeneous setting. This extends
the work \cite{Aetal06} on the reduce max rate policy, where only the existence and uniqueness
of a fixed point was shown but the stability in the deterministic setting was not investigated.
We conclude the paper in Section~\ref{sec:conc} with future research directions.

\section{Model and problem formulation}

Let us consider a Generalized Additive Increase Multiplicative Decrease (G-AIMD) dynamics
with $N$ users in continuous time.
In the absence of control signal, the allocation to user $k$ (e.g., transmission rate
in Internet congestion control or instantaneous power in charging stations for electric
vehicles) increases according to the differential equation:
\begin{equation}
\label{eq:incrdyn}
\frac{dx_k}{dt} = a_k x_k^{\gamma_k},
\end{equation}
with $\gamma_k \in [0,1]$ and $a_k > 0$. Continuous-time models represent well the
TCP sending rate evolution on the scale of several round-trip times \cite{AAP10,Zetal10}.

We consider impulsive control. Namely, when the control signal (impulse)
is sent to user $k$ at time $t$, the resource allocation to user $k$
drastically decreases according to
\begin{equation}
\label{eq:decrdyn}
x_k(t+0) = b_k x_k(t),
\end{equation}
with $b_k \in (0,1)$. We note that the above dynamics is fairly general and covers
at least three important particular cases: if $\gamma_k=0$ we retrieve
the classical Additive Increase Multiplicative Decrease (AIMD) mechanism \cite{S93},
$\gamma_k=3/4$ corresponds to Compound TCP \cite{Tetal06} when queueing delays are not large,
and $\gamma_k=1$ corresponds to the Multiplicative Increase
Multiplicative Decrease (MIMD) mechanism or Scalable TCP \cite{K03}.
MIMD is a very aggressive dynamics \cite{AAP05} and, in contrast,
AIMD is much more gentle. Compound TCP is designed to represent a good balance
between the two extremes.

Let us define formally a class of policies slightly larger than the class of
purely deterministic policies. The need for such a class of policies will be clear from
the subsequent development.

\begin{defi}
\begin{itemize}
\item[(a)]
Let $k=1,2,\dots,N$ be fixed. For user $k,$ a policy is a sequence, say $T=\{T_1,T_2,...\}$ of time moments, when an impulse (a multiplicative decrease in his sending rate) is applied. Here $\{T_i\}$ is a monotone nondecreasing sequence of constants in $[0,\infty]$. It is possible that multiple impulses are applied at a single time moment, but we require $\lim_{n\rightarrow\infty}T_n=\infty$.
\item[(b)] Let $T=\{T_1,T_2,...\}$ and $S=\{S_1,S_2,...\}$ be two policies for user $k$. Then, for each $\beta\in[0,1],$ we denote by $(\beta,T,S)$
a mixture of the two policies, which with probability $\beta$ chooses the sequence $T$ and with the
complementary probability $1-\beta$ chooses the sequence $S$. For user $k$, we denote by
${\cal U}'_k$ (resp., ${\cal U}_k$) the set of policies (resp. all such mixed policies for all $\beta\in[0,1]$).
\item[(c)] We introduce the notation ${\cal U}'=\prod_{k=1}^N{\cal U}'_k$ and ${\cal U}=\prod_{k=1}^N {\cal U}_k .$
\end{itemize}
\end{defi}

Let $k=1,2,\dots,N$ be fixed. Each policy $u_k=(\beta,T,S)\in {\cal U}_k$ defines the dynamics of $x_k(t)$
(stochastic, if $\beta \neq 0,1$),
and the corresponding expectation is denoted by $E^{u_k}[\cdot]$.

Let us denote by $x(t)=[x_1(t) \ \cdots \ x_N(t)]^T$ the vector of resource allocations
at time $t$. Ideally, at each time moment we aim to operate the system under the constraint:
\begin{equation}
\label{eq:capcon}
\sum_{k=1}^{N} x_k(t) \le c, \quad \forall t,
\end{equation}
where $c>0$ is the resource (e.g., transmission capacity or electric power).
It appears that if we substitute the above hard constraint with a soft time-averaged
constraint, the problem becomes more tractable. Namely, consider
$$
\sum_{k=1}^{N} \limsup_{\tau \to \infty} \frac{1}{\tau} E^{u_k}\left[ \int_0^\tau x_k(t) dt \right] \le c.
$$

Our first objective is to propose an impulsive control in closed form, which solves the following constrained problem:
\begin{eqnarray}
\label{eq:alphafair}
&&J(u):= \sum_{k=1}^{N} \liminf_{\tau \to \infty} \frac{1}{\tau} E^{u_k}\left[ \int_0^\tau  \frac{x^{1-\alpha}_k(t)}{1-\alpha} dt \right]
\longrightarrow  \sup_{u\in {\cal U}},\nonumber\\
&&\mbox{subject to:~} \sum_{k=1}^{N} \limsup_{\tau \to \infty} \frac{1}{\tau} E^{u_k}\left[ \int_0^\tau x_k(t) dt \right] \le c.
\end{eqnarray}
Here the initial state is arbitrarily fixed, and will not be indicated, $\alpha>0,$ and the generic notation $u=(u_1,u_2,\dots,u_N)\in {\cal U}$ is in use.
For each $\alpha>0,$ $J(u)$ is the short-term $\alpha$-fairness \cite{AAR12}.
The short-term $\alpha$-fairness is a versatile fairness concept, which retrieves as particular cases:
proportional fairness ($\alpha \to 1$), delay-based fairness ($\alpha=2$) and max-min
fairness ($\alpha \to \infty$).

In order to deal with the control problem under constraints, we use the multiobjective optimization approach.
To this end, let us define the two competing objectives:
\begin{eqnarray*}
&&J(u) = \sum_{k=1}^{N} J_k(u_k)\\
&: =& \sum_{k=1}^{N} \liminf_{\tau \to \infty} \frac{1}{\tau} E^{u_k}\left[ \int_0^\tau  \frac{x^{1-\alpha}_k(t)}{1-\alpha} dt \right]
\ \to \ \sup_{u\in {\cal U}},\\
&&G(u) = \sum_{k=1}^{N} G_k(u_k)\\
&: =& \sum_{k=1}^{N} \limsup_{\tau \to \infty} \frac{1}{\tau} E^{u_k}\left[ \int_0^\tau x_k(t) dt \right]
\quad \to \quad \inf_{u\in {\cal U}}.
\end{eqnarray*}
It appears that it is more convenient to consider $-J$ instead of $J$, which leads to the standard multiobjective problem:
\begin{eqnarray*}
-J(u) \ \to \ \inf_{u\in {\cal U}}, \quad G(u) \ \to \ \inf_{u\in {\cal U}}.
\end{eqnarray*}

Throughout this paper, we assume the following
\begin{assumption}\label{AssumptionA}
$\gamma_k\in[0,1],~ \alpha\ne 1,~ 2-\alpha-\gamma_k\ne 0,~b_k\in(0,1)$ for each $k=1,2,\dots, N.$
\end{assumption}
The particular cases excluded by the assumption can be separately analyzed using similar techniques.
We exclude such cases for the sake of presentation smoothness and because of space limitation.

\section{Control in the relaxed case}
\label{sec:relaxed}

Let us formally justify the reduction of the problem with the relaxed constraint to the
multiobjective formulation and demonstrate how the original solution can be reconstructed.

To scalarize the multiobjective problem, we introduce the variable weight $\lambda \in (0,+\infty)$
and consider the combined criterion
\begin{eqnarray*}
&&L(\lambda,u) =  \sum_{k=1}^{N} L_k(\lambda,u_k)\\
&:=& \sum_{k=1}^{N} (-J_k(u_k)+\lambda G_k(u_k))
\to  \inf_{u\in {\cal U}}.
\end{eqnarray*}

Note that the above problem reduces to $N$ subproblems: for each $k=1,2,\dots,N,$
\begin{eqnarray}\label{NewNn}
-J_k(u_k)+\lambda G_k(u_k)\to \quad \inf_{u_k\in {\cal U}_k}.
\end{eqnarray}

\begin{lemma}\label{Lemma01}
For each $k=1,2,\dots, N,$ and $\lambda>0,$ an optimal policy for problem (\ref{NewNn}) is of threshold type
with the threshold $\bar{x}_k(\lambda)$ given by
\begin{equation}
\label{eq:thres}
\bar{x}_k(\lambda) = \left\{
\frac{(2-\gamma_k)(1-b^{2-\alpha-\gamma_k}_k)}{(1-b^{2-\gamma_k}_k)(2-\alpha-\gamma_k)\lambda}
\right\}^\frac{1}{\alpha}.
\end{equation}
In greater details, under this threshold policy, the user $k$ decreases the sending rate at time $t$ as soon as $x_k(t)\ge \bar{x}_k$. (It is clear that this threshold policy induces a policy in ${\cal U}'_k.$)
\end{lemma}

\par\noindent\textit{Proof.} Let $k=1,2,\dots,N$ be fixed. As was shown in \cite{Aetal15} (see there Theorem~3.1), the policy say $u_k^\ast\in {\cal U}'_k$ defined by the threshold $\bar{x}_k$ given by (\ref{eq:thres}) is optimal to the following problem
\begin{equation}
\label{eq:Lk}
\limsup_{\tau \to \infty} \frac{1}{\tau}
 E^{u_k}\left[ \int_0^\tau \left( -\frac{x^{1-\alpha}_k(t)}{1-\alpha}
+ \lambda x_k(t) \right) dt  \right]
\to \inf_{{\cal U}'_k}.
\end{equation}
It is clear that this policy is also optimal to the above problem but out of all the policies ${\cal U}_k,$ i.e.,
\begin{eqnarray*}
\limsup_{\tau \to \infty} \frac{1}{\tau}
E^{u_k}\left[ \int_0^\tau \left( -\frac{x^{1-\alpha}_k(t)}{1-\alpha}
+ \lambda x_k(t) \right) dt \right]
\to \inf_{{\cal U}_k}.
\end{eqnarray*}
(In fact, if it is outperformed by a mixed policy, then there must be another deterministic policy outperforming this threshold policy, which contradicts the optimality of the threshold policy out of ${\cal U}'_k.$)
Then
\begin{eqnarray*}
&&-J_k(u_k^\ast)+\lambda G_k(u_k^\ast)\\
&=&\lim_{\tau \to \infty} \frac{1}{\tau}
E^{u_k^\ast}\left[ \int_0^\tau \left( -\frac{x^{1-\alpha}_k(t)}{1-\alpha}
+ \lambda x_k(t) \right) dt \right]\\
&\le&\limsup_{\tau \to \infty} \frac{1}{\tau}
E^{u_k}\left[ \int_0^\tau \left( -\frac{x^{1-\alpha}_k(t)}{1-\alpha}
+ \lambda x_k(t) \right) dt \right]\\
&\le&-J_k(u_k)+\lambda G_k(u_k)
\end{eqnarray*}
for each $u_k\in {\cal U}_k$.
$\hfill\Box$

\begin{lemma}
For each $k=1,2,\dots,N,$ and $\lambda>0,$ under the threshold policy $u_k^\ast$ given by (\ref{eq:thres}),  if $\gamma_k<1,$
\begin{eqnarray}\label{Star1}
&&-J_k^*(\lambda):=-J_k(u_k^\ast)\nonumber\\
&=&
-[\bar{x}_k(\lambda)]^{1-\alpha} \frac{(1-b^{2-\alpha-\gamma_k}_k)(1-\gamma_k)}{(1-\alpha)(2-\alpha-\gamma_k)(1-b^{1-\gamma_k}_k)},\nonumber\\
&&G_k^*(\lambda):=G_k(u_k^\ast)\nonumber\\
&=&\bar{x}_k(\lambda) \frac{(1-b^{2-\gamma_k}_k)}{(2-\gamma_k)} \frac{(1-\gamma_k)}{(1-b^{1-\gamma_k}_k)},\nonumber
\\
&&L^*_k(\lambda):=J_k(u_k^\ast)+\lambda G_k(u_k^\ast)\nonumber\\
&=&
- \bar{x}_k(\lambda) \lambda \frac{\alpha}{1-\alpha} \frac{(1-b^{2-\gamma_k}_k)}{(2-\gamma_k)} \frac{(1-\gamma_k)}{(1-b^{1-\gamma_k}_k)};\nonumber\\
\end{eqnarray}
and if $\gamma_k=1,$
\begin{eqnarray}\label{Star2}
&&-J_k^*(\lambda)=
[\bar{x}_k(\lambda)]^{1-\alpha} \frac{1-b_k^{1-\alpha}}{(1-\alpha)^2\ln b_k},\nonumber\\
&&G_k^*(\lambda)=
\bar{x}_k(\lambda) \frac{b_k-1}{\ln b_k},\nonumber
\\
&&L^*_k(\lambda)=
- \bar{x}_k(\lambda) \lambda \frac{\alpha}{1-\alpha} \frac{b_k-1}{\ln b_k}.\nonumber\\
\end{eqnarray}
\end{lemma}

\par\noindent\textit{Proof.} The details needed for the derivation of the above objectives can be found in \cite{Aetal15}.
$\hfill\Box$

Let us now investigate the trade off $G^*_k$ against $-J^*_k$.
We consider two cases (a) $\alpha < 1$ and (b) $\alpha > 1$ separately. The following two observations hold for $\gamma_k\in[0,1].$
\begin{itemize}
\item[(a)] $\alpha < 1$: by equation (\ref{eq:thres}), if $\lambda \to 0$ then $\bar{x}_k \to \infty$ and
consequently $G_k^* \to +\infty$; at the same time, $\lambda^{1-1/\alpha} \to \infty$ and consequently $-J^*_k \to -\infty$.
Now if $\lambda \to \infty$ then $\bar{x}_k \to 0$ and $G_k^* \to 0$; and at the same time
$\lambda^{1-1/\alpha} \to 0$ and consequently $-J^*_k \to 0$.

\item[(b)] $\alpha > 1$: Again by equation (\ref{eq:thres}), if $\lambda \to 0$ then $\bar{x}_k \to \infty$ and
consequently $G_k^* \to +\infty$. However, in this case $\lambda^{1-1/\alpha} \to 0$ and consequently $-J^*_k \to 0$.
Now if $\lambda \to \infty$ then $\bar{x}_k \to 0$ and $G_k^* \to 0$; and at the same time
$\lambda^{1-1/\alpha} \to \infty$ and consequently $-J^*_k \to +\infty$.
\end{itemize}

Next we establish the convexity of the epigraph.

\begin{lemma}
For each $k=1,2,\dots,N,$ $G_{k}^\ast$, legitimately regarded as a function of $-J_k^\ast$, is convex. Moreover,
its epigraph coincides with the convex hull of its graph.
\end{lemma}

\par\noindent{\it Proof:} To prove the convexity, it will be more convenient to consider
the parametrization with respect to $\bar{x}_k$. We note that since there is a one-to-one
correspondence between $\lambda$ and $\bar{x}_k$, the two parametrizations are equivalent.
Observe that
\begin{eqnarray*}
&&-J^*_k(\bar{x}_k) = -c_1 \frac{\bar{x}_k^{1-\alpha}}{1-\alpha}, \quad c_1>0,\\
&&G^*_k(\bar{x}_k) = c_2 \bar{x}_k, \quad c_2>0,
\end{eqnarray*}
where the constants $c_1,c_2$ come from (\ref{Star1}) (resp., (\ref{Star2})) when $\gamma_k<1$ (resp. $\gamma_k=1$).
Thus, we can write
\begin{eqnarray*}
G^*_k(-J^*_k) = c_2 \left[\frac{(-J^*_k)(1-\alpha)}{-c_1}\right]^{\frac{1}{1-\alpha}}.
\end{eqnarray*}
Hence,
\begin{eqnarray*}
\frac{dG^*_k}{d(-J^*_k)} =
- \frac{c_2}{c_1} \left[\frac{(-J^*_k)(1-\alpha)}{-c_1}\right]^{\frac{\alpha}{1-\alpha}},
\end{eqnarray*}
and
\begin{eqnarray*}
\frac{d^2G^*_k}{d(-J^*_k)^2} =
\alpha \frac{c_2}{c_1^2} \left[\frac{(-J^*_k)(1-\alpha)}{-c_1}\right]^{\frac{2\alpha-1}{1-\alpha}} > 0,
\end{eqnarray*}
since $[(-J^*_k)(1-\alpha)/(-c_1)]>0$ always.

\begin{figure}[h]
\begin{center}
\subfigure[Case $\alpha < 1$.]{\includegraphics[scale=0.28]{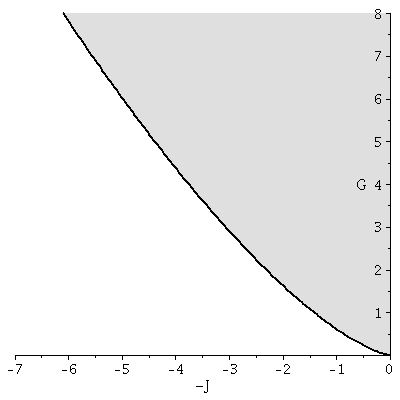}}
\subfigure[Case $\alpha > 1$.]{\includegraphics[scale=0.28]{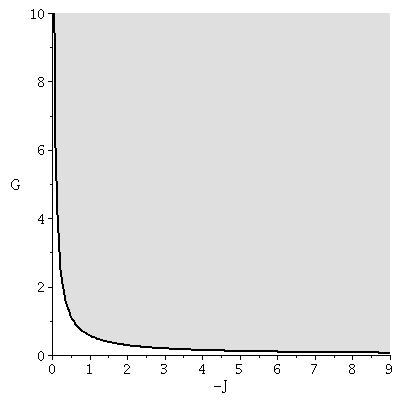}}\\
\caption{Epigraph of $\Omega_k$.}
\vspace{-0.3cm}
\label{fig:Epigraph}
\end{center}
\end{figure}

The last assertion follows from the two observations before this lemma: there is no asymptote if $\alpha<1$, and the same conclusion holds if $\alpha>1$, too. Examples of the epigraph in the two cases  are displayed in Figures~\ref{fig:Epigraph}.(a)~and~\ref{fig:Epigraph}.(b). This completes the proof.
\hfill $\Box$

%

\begin{remark}\label{Remark01}
For each $k=1,2,\dots,N,$ denote by $\Omega_k$ the convex hull of the graph (or equivalently the epigraph, according to the previous lemma) of $G^*_k$ as a function of $-J^*_k$. It can be seen that $\{(-J_k(u_k),G_k(u_k)):u_k\in {\cal U}_k\}=\Omega_k$. Indeed, if there is some $u_k\in {\cal U}_k$ such that $((-J_k(u_k),G_k(u_k)))\notin \Omega_k$, it can only lie below the graph of $G^\ast_k$ against $-J_k^\ast$, but then for some $\lambda>0,$ it contradicts the fact that the threshold policy $u_k^\ast$ given by $\bar{x}_k(\lambda)$ is optimal for problem (\ref{NewNn}) with the same $\lambda>0.$ This observation is important for the argument below.
\end{remark}

Now we consider problem (\ref{eq:alphafair}), and reformulate it in the space of performance vectors. That is, we reformulate
\begin{eqnarray*}
\left\{
\begin{array}{l}
-J(u) = \sum_{k=1}^{N} (-J_k(u_k)) \to \inf_{u \in {\cal U}},\\
G(u) - c = \sum_{k=1}^{N} G_k(u_k) - c \le 0,
\end{array}
\right.
\end{eqnarray*}
as
\begin{eqnarray}\label{PerformanceVectorProb}
\left\{
\begin{array}{l}
- \tilde J(\omega) \to \inf_{\omega \in \Omega},\\
\tilde G(\omega) - c \le 0,
\end{array}
\right.
\end{eqnarray}
where
\begin{eqnarray*}
&&- \tilde J(\omega): = \sum_{k=1}^N \omega^1_k,\\
&&\tilde G(\omega) = \sum_{k=1}^N \omega^2_k,
\end{eqnarray*}
and where
\begin{eqnarray*}
\omega = \{(\omega^1_k,\omega^2_k)\}_{k=1}^N \in \Omega: = \prod_{k=1}^N \Omega_k \subset \mathbb{R}^{2N}.
\end{eqnarray*}
In fact, these two problems are equivalent because of the following. For each $u \in {\cal U},$ there exists some $\omega \in \Omega$ such that $J(u)=\tilde J(\omega)$ and $G(u)=\tilde G(\omega)$ , and
conversely for each $\omega \in \Omega,$ there exists some $u \in {\cal U}$ satisfying  $\tilde J(\omega)=J(u)$ and $\tilde G(\omega)=G(u)$; recall Remark \ref{Remark01}. However,
the correspondence $u \leftrightarrow \omega$ may be not one-to-one.

We shall effectively solve problem (\ref{PerformanceVectorProb}), whose optimal solution then induces one to problem (\ref{eq:alphafair}).

The main statement is now in position.
\begin{theorem}\label{Theorem001}
The following assertions hold.
\begin{itemize}
\item[(a)] The set $\Omega$ is convex in $\mathbb{R}^{2N}$, the functions $-\tilde{J}(\omega)$ and $\tilde{G}(\omega)$ on $\Omega$ are convex and real-valued.
\item[(b)] There exists some $\hat{\omega}\in \Omega$ such that $\tilde{G}(\hat{\omega})<c$, i.e., Slater's condition for problem (\ref{PerformanceVectorProb}) is satisfied.
\item[(c)] The threshold policy $u^\ast=(u^\ast_1,\dots,u^\ast_N)\in {\cal U}'$ is optimal for problem (\ref{eq:alphafair}), where for each $k=1,\dots,N,$ $u_k^\ast\in {\cal U}'_k$ is induced by the threshold $\bar{x}_k(\lambda^\ast)$, with
\begin{eqnarray}\label{eq:mult_nonsym}
&&\lambda^*= \frac{1}{c^{\alpha}}
\biggl(\sum_{k=1,\dots, N:\gamma_k\ne 1}\frac{(1-\gamma_k)}{(2-\gamma_k)} \frac{(1-b^{2-\gamma_k}_k)}{(1-b^{1-\gamma_k}_k)}\nonumber\\
&&\times\left(\frac{2-\gamma_k}{1-b^{2-\gamma_k}_k}\right)^{1/\alpha}
\left(\frac{1-b^{2-\alpha-\gamma_k}_k}{2-\alpha-\gamma_k}\right)^{1/\alpha}\nonumber\\
&&+\sum_{k=1,\dots,N:\gamma_k=1} \left(\frac{1-b_k^{1-\alpha}}{1-\alpha}\right)^{1/\alpha}
\frac{(1-b_k)^{(1-\alpha)/\alpha}}{(-\ln b_k)}
\biggr)^{\alpha}\nonumber\\
\end{eqnarray}

In the homogeneous case ($\gamma_k=\gamma$ and $b_k=b$) the expression becomes even simpler:
\begin{eqnarray*}
\lambda^*
&=&\frac{N^{\alpha}}{c^{\alpha}}
\frac{(1-\gamma)^\alpha}{(1-b^{1-\gamma})^\alpha}
\frac{(1-b^{2-\gamma})^\alpha}{(2-\gamma)^\alpha}\\
&&
\times\frac{(2-\gamma)}{(1-b^{2-\gamma})}\frac{(1-b^{2-\alpha-\gamma})}{(2-\alpha-\gamma)},
\end{eqnarray*}
and, consequently,
\begin{equation}
\label{eq:xbarstar}
\bar{x}_k(\lambda^*) = \frac{c}{N} \frac{1-b^{1-\gamma}}{1-\gamma} \frac{2-\gamma}{1-b^{2-\gamma}},
\end{equation}
in case $\gamma \ne 1,$ and
\begin{eqnarray*}
\lambda^\ast&=& \frac{1}{c^\alpha} \left( N \frac{(1-b^{1-\alpha})^{1/\alpha}(b-1)}{((1-b)(1-\alpha))^{1/\alpha} \ln b}  \right)^\alpha.
\end{eqnarray*}
in case $\gamma =1.$
\end{itemize}

\end{theorem}
\par\noindent\textit{Proof.}
 Part (a) is evident. For part (b), note that one can take such $\hat\omega\in \Omega$ that
\begin{eqnarray*}
\forall k \in \{1,...,N\}, \quad \hat\omega^2_k < \frac{c}{N}.
\end{eqnarray*}
This is possible because $G^*_k(\lambda)$ approaches zero when $\lambda \to \infty$.
Thus, Slater's condition is satisfied.

The rest of this proof verifies part (c). For each $\lambda>0,$ let $\omega^\ast(\lambda)\in \Omega$ be generated by the threshold policy $u=(u_1,\dots,u_N)$ determined by the threshold $\bar{x}_k(\lambda)$, $k=1,\dots,N.$ We solve
\begin{eqnarray}\label{002}
\tilde{G}(\omega^\ast(\lambda))=c
\end{eqnarray}
for $\lambda^\ast$ given by (\ref{eq:mult_nonsym}).
Then it holds that
\begin{eqnarray}\label{001}
&&-\tilde{J}(\omega^\ast(\lambda^\ast))+\lambda^\ast (\tilde{G}(\omega^\ast(\lambda^\ast))-c)\nonumber\\
&\le&
-\tilde{J}(\omega)+\lambda^\ast (\tilde{G}(\omega)-c),~\forall~\omega\in \Omega,
\end{eqnarray}
by Lemma \ref{Lemma01}. According to Theorem 1 of Section 8.4 in \cite{L69}, this shows that
$\omega^\ast(\lambda^\ast)\in \Omega$ solves problem (\ref{PerformanceVectorProb}). Part (c) immediately follows. $\hfill\Box$

Consider $\lambda^\ast$ given by (\ref{eq:mult_nonsym}). According to (\ref{001}) and that (\ref{002}) is satisfied by $\omega^\ast(\lambda^\ast)$, we see
\begin{eqnarray*}
\mu_0=\inf_{\omega \in \Omega} \{-\tilde{J}(\omega) + \lambda^* (\tilde{G}(\omega)-c)\},
\end{eqnarray*}
where
\begin{eqnarray*}
\mu_0 = \inf -\tilde{J}(\omega), \ \mbox{subject to} \ \omega \in \Omega, \ \tilde{G}(\omega) \le c.
\end{eqnarray*}
Any constant $\lambda\ge 0$ satisfying the above equality with $\lambda^\ast$ being replaced by $\lambda$ is sometimes called a geometric multiplier for problem (\ref{PerformanceVectorProb}), see Definition 6.1.1 of \cite{Bert}. The following result from \cite[Thm.~1 in Sect.~8.3]{L69}, see the proof therein, shows that $\lambda^\ast$ is the unique geometric multiplier for problem (\ref{PerformanceVectorProb}).

\begin{prop}
\label{thm:Luen}
Let $\Omega \in \mathbb{R}^m$ be a convex set. Let $f$ be a real-valued convex function on $\Omega$
and $G$ be a real-valued convex function on $\Omega$. Assume the existence of a point
$\hat\omega \in \Omega$ for which $G(\hat\omega) < 0$. Let
\begin{equation}
\label{eq:thm1.1}
\mu_0 = \inf f(\omega), \ \mbox{subject to} \ \omega \in \Omega, \ G(\omega) \le 0,
\end{equation}
and assume $\mu_0$ is finite. Then there is a number $\lambda' \ge 0$ such that
\begin{equation}
\label{eq:thm1.2}
\mu_0 = \inf_{\omega \in \Omega} \{f(\omega) + \lambda' G(\omega)\},
\end{equation}
and thus a geometric multiplier exists.
Furthermore, for each geometric multiplier $\lambda',$ if the infimum is achieved in (\ref{eq:thm1.1}) by an $\omega^* \in \Omega$,
$G(\omega^*) \le 0$, it is achieved by $\omega^*$ in (\ref{eq:thm1.2}) and
\begin{eqnarray*}
\lambda' G(\omega^*) = 0.
\end{eqnarray*}
\end{prop}

\begin{corollary}
$\lambda^\ast$ given by (\ref{eq:mult_nonsym}) is the unique geometric multiplier for problem (\ref{PerformanceVectorProb})
\end{corollary}
\par\noindent\textit{Proof.} Suppose $\lambda'\ge 0$ is a geometric multiplier for problem (\ref{PerformanceVectorProb}). Let $\omega^\ast=\omega^\ast(\lambda^\ast)\in \Omega$ be as in the proof of Theorem \ref{Theorem001}. Let us verify that $\lambda'>0.$
Suppose for contradiction that $\lambda'=0$.
Consider the case of $\alpha < 1.$ Remember, $\mu_0$ is finite. However, since $\lambda'=0$ is a geometric multiplier,
\begin{eqnarray*}
\forall k \in \{1,...,N\}, \quad \inf_{\omega_k \in \Omega_k} [\omega^1_k] = -\infty,
\end{eqnarray*}
and
\begin{eqnarray*}
\mu_0 = \inf_{\omega \in \Omega} [-\tilde J(\omega)] = -\infty,
\end{eqnarray*}
which leads to a contradiction. Consider the case of $\alpha > 1$. Then $\mu_0$ is strictly positive.
However,
\begin{eqnarray*}
\forall k \in \{1,...,N\}, \quad \inf_{\omega_k \in \Omega_k} [\omega^1_k] = 0,
\end{eqnarray*}
and since $\lambda'=0$ is a geometric multiplier,
\begin{eqnarray*}
\mu_0 = \inf_{\omega \in \Omega} [-\tilde J(\omega)] = 0,
\end{eqnarray*}
which again leads to a contradiction.

Thus, $\lambda'>0.$ According to Proposition \ref{thm:Luen}, $\lambda'$ satisfies $\tilde{G}(\omega^\ast(\lambda'))=c$, which admits the unique solution $\lambda'=\lambda^\ast.$  $\ \hfill\Box$

The above results call for a number of distributed control algorithms.
At first, let us suppose that the numbers of users' types are known to all users
or broadcasted to the users by a central authority (e.g., SDN controller).
Then, each user can calculate its threshold $\bar{x}_k$ by (\ref{eq:thres}),(\ref{eq:mult_nonsym})
and can control his rate $x_k(t)$ by reducing it when the threshold is achieved.
Thus, except for the complete initial knowledge of the system's parameters,
no further exchange of information is required.

Then, another interesting case is when each user knows its individual parameters but not the parameters
of the other users. In this case, the central controller can calculate the Lagrange multiplier by
equation (\ref{eq:mult_nonsym}) and distribute it to the users.

\section{Index policy for hard constraint}
\label{sec:hard}

Since $\lambda=\lambda(\bar{x})$ is monotone and decreasing function of $\bar{x}$,
the comparison of $\lambda(\bar{x}_k)$ with $\lambda^*$ provides the optimal solution for the relaxed problem
formulation. What is more, the fact that $\lambda(\cdot)$ is a monotone and decreasing function implies
indexability of the problem with hard constraint \cite{W88}.

Then, we can propose the following heuristic for the case of hard constraint \cite{W88}: whenever
the hard constraint (\ref{eq:capcon}) is achieved, the user with the minimal value of $\lambda(x_k(t))$
reduces his rate. Let us call the resulting policy the Whittle-type index policy or briefly the index policy.

It is very intriguing to observe that
the expression for $\lambda(\bar{x})$ contains neither the parameters of the other users nor the number of users.
Therefore, may be the Whittle index type approach can be very useful in the adaptive scenario when
the number of users changes with time.

From now on, in this section, we consider the homogeneous case, i.e., we suppose $a_k=a,~\gamma_k=\gamma\in[0,1)$ and $\beta_k=b\in(0,1)$ for each $k=1,\dots,N$. This is the standard first step in the analysis of index policies \cite{WW90}.
As previously, Assumption~\ref{AssumptionA} is supposed to hold without explicit references.
It is without loss of generality to assume $a=1$.

Let $u_{ind}$ be the index policy. Let $u^\ast$ be the threshold policy obtained in Theorem \ref{Theorem001}, which is optimal for problem (\ref{eq:alphafair}). Note that the index policy $u_{ind}$ satisfies the hard capacity constraint (\ref{eq:capcon}). Therefore, denoting ${\cal U}_H$ as the class of policies satisfying the hard capacity constraint (\ref{eq:capcon}), one has
\begin{eqnarray*}
J(u_{ind},x,c,N)\le \sup_{u\in {\cal U}_H}J(u,x,c,N)\le J(u^\ast,x,c,N),
\end{eqnarray*}
for each initial state $x$, capacity constraint $c$, and the number of users $N$, which we signify in this section for the following reason. Our objective is to show that the index policy is asymptotically optimal in the following sense:
\begin{eqnarray}\label{AsymNeeded06}
\lim_{N\rightarrow \infty}\frac{1}{N}J(u_{ind},x,cN,N)
=\lim_{N\rightarrow \infty} \frac{1}{N}J(u^\ast,x,cN,N).
\end{eqnarray}
In the important case of $\gamma=0$ corresponding to the AIMD dynamics, we show that the index policy $u_{ind}$ is asymptotically optimal for each initial state, and in case of $\gamma\in[0,1)$, we show that it is asymptotically optimal
for the initial states close enough to the steady state.

\subsection{The AIMD $(\gamma=0)$ case}

Suppose $\gamma_k=0$ and $\beta_k=b\in(0,1)$ for each $k=1,\dots,N$.

We first observe that since $\lambda(\cdot)$ is monotone and decreasing in the homogeneous case
the index policy is equivalent to the policy that reduces the maximal sending rate at the moment
when the hard constraint is achieved.
Let us now consider, under the index policy, the sequence of the sending rates, observed at each time
when the capacity constraint is met. Following \cite{Aetal06},
for each $\tilde{x}=(\tilde{x}_1,\dots,\tilde{x}_N)\in(0,\infty)^N$  such that
\begin{eqnarray}\label{Needed02}
&&\tilde{x}\ge \tilde{x}_2\ge\dots\ge \tilde{x}_N>0;~\sum_{i=1}^N \tilde{x}_i=c,
\end{eqnarray}
we introduce
\begin{eqnarray*}
g(\tilde{x}):=(g_1(\tilde{x}),\dots,g_N(\tilde{x}))
\end{eqnarray*}
defined in the following way. If
\begin{eqnarray}\label{EqnNeeded}
\tilde{x}_k\ge  b \tilde{x}_1>\tilde{x}_{k+1}
\end{eqnarray}
for some $1\le k\le N$ with the convention $x_{N+1}:=0,$
then
\begin{eqnarray*}
\left\{
\begin{array}{l}
g_1(\tilde{x}):=\tilde{x}_2+b_1 \tilde{x}_1,\\
\vdots\\
g_{k-1}(\tilde{x}):=\tilde{x}_k+b_1\tilde{x}_1,\\
g_{k}(\tilde{x}):=b\tilde{x}_1+b_1\tilde{x}_1,\\
g_{k+1}(\tilde{x}):=\tilde{x}_{k+1}+b_1\tilde{x}_1,\\
\vdots\\
g_N(\tilde{x}):=\tilde{x}_N+b_1\tilde{x}_1,
\end{array}
\right.
\end{eqnarray*}
where the last two lines are not relevant if $k=N,$ and $b_1=\frac{1-b}{N}>0$ is a constant. Note that if $\Delta(\tilde{x})$ denotes the time duration since the reduction of $\tilde{x}$ according to the index policy until the next time when the hard capacity constraint is met, then  \begin{eqnarray}\label{Needed05}
b_1\tilde{x}_1=\frac{(1-b)\tilde{x}_1}{N}=a\Delta(\tilde{x}).
\end{eqnarray}
The interpretation of $g(\tilde{x})$ is the vector of the ordered sending rates from the largest to the smallest one, when the next time the hard capacity constraint is met (before the reduction), starting from $\tilde{x}.$  Put $\tilde{x}=:\tilde{x}^{(0)}$, with $\tilde{x}\in(0,\infty)^N$ is a fixed vector satisfying (\ref{Needed02}),
\begin{eqnarray*}
\tilde{x}^{(m)}:=g(\tilde{x}^{(m-1)})=:g^{(m)}(\tilde{x}),~m\ge 1,
\end{eqnarray*}
we are interested in $\tilde{x}^{(m)}$ as $m\rightarrow \infty.$ Let us introduce for each vector $\tilde{x}$ satisfying (\ref{Needed02})
$k(\tilde{x})$ as the integer $k$ satisfying (\ref{EqnNeeded}).

\begin{theorem}\label{KostiasTheorem}
Suppose $\gamma_k=0$ and $\beta_k=\beta$ for each $k=1,\dots,N$. Then the mapping $g$ has a unique fixed point, say $\tilde{x}^\ast=(\tilde{x}^\ast_1,\dots,\tilde{x}^\ast_N)$, in the space of vectors satisfying (\ref{Needed02}), given by
\begin{equation}
\label{eq:AIMDfixedpoint}
\tilde{x}^\ast_n=\left(b+\frac{(N-n+1)(1-b)}N\right)\frac{c}{Nb+\frac{(N+1)(1-b)}{2}},
\end{equation}
$1\le n\le N$, and $\tilde{x}^{(m)}\rightarrow \tilde{x}^\ast$ as $m\rightarrow \infty.$
\end{theorem}
\par\noindent\textit{Proof.}
Firstly, note that there exists some integer $m_0>0$ such that $k(\tilde{x}^{(m_0)})=N,$ for otherwise, the sending rate of some user would have blown up to $\infty$, violating the hard capacity constraint.

Next, observe that if for some $m\ge 0$,\ $k(\tilde{x}^{(m)})=N,$ then $k(\tilde{x}^{(m+1)})=N$ as well.
Indeed, if this was not the case, then we would have
\begin{eqnarray*}
b \tilde{x}^{(m+1)}_1=b(\tilde{x}_2^{(m)}+b_1\tilde{x}_1^{(m)})>b\tilde{x}_1^{(m)}+b_1\tilde{x}_1^{(m)}
\end{eqnarray*}
and thus
\begin{eqnarray*}
0\ge b(\tilde{x}_2^{(m)}-\tilde{x}_1^{(m)})>b_1(\tilde{x}_1^{(m)}-b\tilde{x}_1^{(m)})>0,
\end{eqnarray*}
which is a desired contradiction. Therefore, for some and all subsequent steps, the maximal sending rate (before reduction) when the hard capacity constraint is met will become the minimal sending rate (just after reduction).

Consequently, for all large enough $m,$ we have
\begin{eqnarray*}
&&\tilde{x}^{(m+1)}=g(\tilde{x}^{(m)})\\
&=&\left(
     \begin{array}{ccccc}
       b_1 & 1 & 0 & \dots & 0 \\
       b_1 & 0 & 1 & \dots & 0 \\
       \vdots & \vdots & \vdots & \dots & \vdots \\
       b_1 & 0 & 0 & \dots & 1 \\
       b_1+b & 0 & 0 & \dots & 0 \\
     \end{array}
   \right) \tilde{x}^{(m)},
\end{eqnarray*}
Since the matrix
\begin{eqnarray*}
A = \left(
     \begin{array}{ccccc}
       b_1 & 1 & 0 & \dots & 0 \\
       b_1 & 0 & 1 & \dots & 0 \\
       \vdots & \vdots & \vdots & \dots & \vdots \\
       b_1 & 0 & 0 & \dots & 1 \\
       b_1+b & 0 & 0 & \dots & 0 \\
     \end{array}
   \right)
\end{eqnarray*}
is an aperiodic irreducible (column) stochastic matrix, we conclude that $\tilde{x}^{(m)}$ converges to the unique fixed point $\tilde{x}^\ast$ of $g$ in the space of vectors satisfying (\ref{Needed02}).

Let us compute the fixed point $\tilde{x}^\ast=(\tilde{x}^\ast_1,\dots,\tilde{x}^\ast_N)$ by solving the following system:
\begin{eqnarray*}
\left\{
\begin{array}{l}
\tilde{x}_1^\ast:=\tilde{x}^\ast_2+b_1 \tilde{x}^\ast_1,\\
\vdots\\
\tilde{x}_2^\ast:=\tilde{x}^\ast_3+b_1\tilde{x}^\ast_1,\\
\vdots\\
\tilde{x}^\ast_N:=b\tilde{x}^\ast_1+b_1\tilde{x}^\ast_1,\\
\sum_{i=1}^N \tilde{x}^\ast_i=c,
\end{array}
\right.
\end{eqnarray*}
which gives
\begin{eqnarray*}
\left\{
\begin{array}{l}
\tilde{x}^\ast_n:=(b+(N-n+1)b_1)\tilde{x}^\ast_1,~1\le n\le N,\\
\sum_{i=1}^N \tilde{x}^\ast_i=c,
\end{array}
\right.
\end{eqnarray*}
Therefore,
\begin{eqnarray*}
&&\tilde{x}^\ast_1=\frac{c}{\sum_{i=1}^N (b+ib_1)}=\frac{c}{Nb+\frac{(N+1)(1-b)}{2}},\\
&&\tilde{x}^\ast_n=(b+\frac{(N-n+1)(1-b)}N)\frac{c}{Nb+\frac{(N+1)(1-b)}{2}},\\
&&~2\le n\le N,
\end{eqnarray*}
see (\ref{eq:AIMDfixedpoint}).
$\hfill\Box$

We remark that the reduce maximal sending rate policy was investigated in \cite{Aetal06}.
There only the existence and uniqueness of the fixed point (\ref{eq:AIMDfixedpoint}) was
established but the convergence or the absence of cycling behaviour was not shown.

Next, we shall scale the capacity constraint $c$ by a multiplicative constant $N$.
When we do such scaling, it is convenient to signify the dependence of $\tilde{x}^\ast(cN)$
on the capacity constraint $cN$ explicitly.

\begin{theorem}
Suppose $\gamma_k=0$ and $\beta_k=\beta$ for each $k=1,\dots,N$. Then for each initial state, the index policy is asymptotically optimal for the problem with hard capacity constraint, i.e.,
(\ref{AsymNeeded06}) holds.
\end{theorem}
\par\noindent\textit{Proof.}  Note that for each fixed $N,$ according to Theorem \ref{KostiasTheorem},
\begin{eqnarray}\label{KKk}
\frac{1}{N}J(u_{ind},x,CN,N)=\frac{\int_{b\tilde{x}^\ast_1(cN)}^{\tilde{x}^\ast_1(cN)}   \frac{x^{1-\alpha}}{1-\alpha}dx }{\int_{b\tilde{x}^\ast_1(cN)}^{\tilde{x}^\ast_1(cN)} 1 dx}.
\end{eqnarray}
Therefore, the statement would be proved if we can verify that as $N\rightarrow\infty,$ $\tilde{x}^\ast_1(cN)\rightarrow  \bar{x}_1(\lambda^\ast)$, with $\bar{x}_1(\lambda^\ast)$ being given by (\ref{eq:thres}), where
\begin{eqnarray*}
\lambda^*
&=&\frac{N^{\alpha}}{(cN)^{\alpha}}
\frac{1}{(1-b)^\alpha}
\frac{(1-b^{2})^\alpha}{2^\alpha}\\
&&
\times\frac{2}{(1-b^{2})}\frac{(1-b^{2-\alpha})}{(2-\alpha)}\\
&=&\frac{1}{c^\alpha}\frac{(1+b)^\alpha}{2^{\alpha-1}}\frac{1-b^{2-\alpha}}{(1-b^2)(2-\alpha)}.
\end{eqnarray*}
is provided by Theorem \ref{Theorem001}(c). This can be done through direct calculations, which we omit, as in the proof of Theorem \ref{KostiaTheorem02} below, this fact will be verified in a more general setup.
$\hfill\Box$

\subsection{The G-AIMD $(\gamma\in[0,1))$ case}
In this subsection, we assume $\gamma_k=\gamma\in[0,1)$ for each $k=1,\dots,N$.

As in the previous subsection, whose notations are adopted here without repeating,
we consider, under the index policy, the sequence of the sending rates, observed at each time when the hard capacity constraint is met (right before the reduction in the maximal sending rate). Now for each $\tilde{x}\in(0,\infty)^N$ satisfying (\ref{Needed02}), if (\ref{EqnNeeded}) holds for some $1\le k\le N$ with the convention $x_{N+1}:=0,$
then $g(\tilde{x})$ is defined by
\begin{eqnarray*}
\left\{
\begin{array}{l}
g_1(\tilde{x}):=[\tilde{x}_2^{1-\gamma}+(1-\gamma)\Delta(\tilde{x})]^{\frac{1}{1-\gamma}},\\
\vdots\\
g_{k-1}(\tilde{x}):=[\tilde{x}_k^{1-\gamma}+(1-\gamma)\Delta(\tilde{x})]^{\frac{1}{1-\gamma}},\\
g_{k}(\tilde{x}):=[(b\tilde{x}_1)^{1-\gamma}+(1-\gamma)\Delta(\tilde{x})]^{\frac{1}{1-\gamma}},\\
g_{k+1}(\tilde{x}):=[\tilde{x}_{k+1}^{1-\gamma}+(1-\gamma)\Delta(\tilde{x})]^{\frac{1}{1-\gamma}},\\
\vdots\\
g_N(\tilde{x}):=[\tilde{x}_N^{1-\gamma}+(1-\gamma)\Delta(\tilde{x})]^{\frac{1}{1-\gamma}},
\end{array}
\right.
\end{eqnarray*}
where the last two lines are not relevant if $k=N,$ and $\Delta(\tilde{x})$ is such that
$\sum_{k=0}^{N} g_k(\tilde{x}) = c$.

As in the proof of Theorem \ref{KostiasTheorem}, there exists some integer $m_0>0$ such that $k(\tilde{x}^{(m_0)})=N.$ Moreover,
if for some $m\ge 0$, $k(\tilde{x}^{(m)})=N,$ then $k(\tilde{x}^{(m+1)})=N$ as well, for otherwise, we would have
\begin{eqnarray*}
&&b \tilde{x}_1^{(m)}=b[(\tilde{x}_2^{(m)}){1-\gamma}+(1-\gamma)\Delta(\tilde{x}^{(m)})]^{\frac{1}{1-\gamma}}\\
&>&[(b\tilde{x}_1^{(m)}){1-\gamma}+(1-\gamma)\Delta(\tilde{x}^{(m)})]^{\frac{1}{1-\gamma}}
\end{eqnarray*}
so that
\begin{eqnarray*}
&&0\ge b^{1-\gamma}[(\tilde{x}_2^{(m)})^{1-\gamma}-(\tilde{x}_1^{(m)})^{1-\gamma}]\\
&>&(1-b^{1-\gamma})(1-\gamma)\Delta(\tilde{x}^{(m)})>0,
\end{eqnarray*}
which is a contradiction. Therefore, for all large enough $m$,
\begin{equation}
\label{eq:nonlindyn}
\left\{
\begin{array}{l}
g_1(\tilde{x}^{(m+1)}):=[(\tilde{x}_2^{(m)})^{1-\gamma}+(1-\gamma)\Delta(\tilde{x}^{(m)})]^{\frac{1}{1-\gamma}},\\
\vdots \\
g_{k}(\tilde{x}^{(m+1)}):=[(\tilde{x}_{k+1}^{(m)})^{1-\gamma}+(1-\gamma)\Delta(\tilde{x}^{(m)})]^{\frac{1}{1-\gamma}},\\
\vdots \\
g_N(\tilde{x}^{(m+1)}):=[(b\tilde{x}_1^{(m)})^{1-\gamma}+(1-\gamma)\Delta(\tilde{x}^{(m)})]^{\frac{1}{1-\gamma}}.
\end{array}
\right.
\end{equation}

Next note that the unique fixed point $\tilde{x}^\ast$ of $g$ in the space of vectors satisfying (\ref{Needed02})
can be computed by solving
\begin{eqnarray*}
\left\{
\begin{array}{l}
\tilde{x}^\ast_1:=[(\tilde{x}_2^{\ast})^{1-\gamma}+(1-\gamma)\Delta(\tilde{x}^{\ast})]^{\frac{1}{1-\gamma}},\\
\vdots\\
\tilde{x}^\ast_k:=[(\tilde{x}_{k+1}^{\ast})^{1-\gamma}+(1-\gamma)\Delta(\tilde{x}^{\ast})]^{\frac{1}{1-\gamma}},\\
\vdots\\
\tilde{x}^\ast_N:=[(b\tilde{x}_1^{\ast})^{1-\gamma}+(1-\gamma)\Delta(\tilde{x}^{\ast})]^{\frac{1}{1-\gamma}}.\\
\sum_{i=1}^N \tilde{x}_i^\ast=c,
\end{array}
\right.
\end{eqnarray*}
and hence, we have
\begin{eqnarray}
\label{eq:nonlinfixedpoint}
&&\tilde{x}^\ast_1=\frac{c}{\sum_{i=1}^N (b^{1-\gamma} +\frac{i}{N}(1-b^{1-\gamma}) )^{\frac{1}{1-\gamma}}},\\
&&\tilde{x}^\ast_n=\tilde{x}^\ast_1 [b^{1-\gamma}+\frac{N-n+1}{N}(1-b^{1-\gamma})]^{\frac{1}{1-\gamma}}, \nonumber\\
&&~2\le n\le N. \nonumber
\end{eqnarray}

Let us now show the local asymptotic stability of the nonlinear dynamics (\ref{eq:nonlindyn}).
After linearization around the fixed point (\ref{eq:nonlinfixedpoint}) and change of variables
$\delta \tilde{z}_k = \delta \tilde{x}_k/(\tilde{x}^\ast_k)^\gamma$, $\delta \tilde{x}_k=\tilde{x}_k-\tilde{x}^\ast_k$,
the linearized dynamics (\ref{eq:nonlindyn}) takes the form
$$
\delta \tilde{z}^{(m+1)} = [I-\underline{1}p^T]B \delta \tilde{z}^{(m)},
$$
where
$p_k=(\tilde{x}^\ast_k)^\gamma/((\tilde{x}^\ast_1)^\gamma+...+(\tilde{x}^\ast_N)^\gamma)$, for $k=1,...,N,$\\
$\underline{1}$ is the vector of ones, and the matrix $B$ is defined as follows:
\begin{eqnarray*}
B = \left(
     \begin{array}{ccccc}
       0 & 1 & 0 & \dots & 0 \\
       0 & 0 & 1 & \dots & 0 \\
       \vdots & \vdots & \vdots & \dots & \vdots \\
       0 & 0 & 0 & \dots & 1 \\
       b^{1-\gamma} & 0 & 0 & \dots & 0 \\
     \end{array}
   \right).
\end{eqnarray*}
Let us investigate the eigenvalues of $[I-\underline{1}p^T]B$. First, we note that it has
one zero eigenvalue with the associated left eigenvector $p^T$. Then, let us multiply
$x^T[I-\underline{1}p^T]B=\lambda x^T$ from the right by the right null eigenvector
$[b^{-(1-\gamma)} \ 1 \ \cdots \ 1]^T$, which gives
$$
\lambda [b^{-(1-\gamma)} \ 1 \ \cdots \ 1]^T x = 0.
$$
Thus, if $[b^{-(1-\gamma)} \ 1 \ \cdots \ 1]^T x \neq 0$, $\lambda=0$ and on the other hand,
if $\lambda \neq 0$, we have $[b^{-(1-\gamma)} \ 1 \ \cdots \ 1]^T x =0$.

Using the additional condition $[b^{-(1-\gamma)} \ 1 \ \cdots \ 1]^T x =0$, after some
algebra, we obtain an equation for the nonzero eigenvalues of $[I-\underline{1}p^T]B$:
$$
\lambda^N-b^{1-\gamma} + (1-b^{1-\gamma})(p_1+p_2\lambda+...+p_n\lambda^{n-1}) = 0,
$$
which is in fact also the characteristic equation of a column stochastic Leslie matrix.
Consequently, $[I-\underline{1}p^T]B$ has $N-1$ roots smaller than unity and one zero root,
which implies local asymptotic stability.

If we now replace $c$ with $c_1 N$ in the expression (\ref{eq:nonlinfixedpoint}) for $\tilde{x}_1^\ast$,
then
\begin{eqnarray*}
&&\lim_{N\rightarrow \infty}\tilde{x}^\ast_1(c_1 N)=\frac{c_1}{\frac{1}{N}\sum_{i=1}^N (b^{1-\gamma} +\frac{i}{N}(1-b^{1-\gamma}) )^{\frac{1}{1-\gamma}}}\\
&=& \frac{c_1}{\int_0^1 (b^{1-\gamma}+x(1-b^{1-\gamma}))^{\frac{1}{1-\gamma}}dx}\\
&=&c_1\frac{1-b^{1-\gamma}}{1-\gamma}\frac{2-\gamma}{1-b^{2-\gamma}}=\bar{x}_1(\lambda^\ast(cN)),
\end{eqnarray*}
which is the expression (\ref{eq:xbarstar}) with $c$ being replaced by $c_1 N$.
Thus, we have the following result:

\begin{theorem}\label{KostiaTheorem02}
Suppose $\gamma_k=\gamma\in[0,1)$ for each $k=1,\dots,N$. Then, if the initial state is close enough to the steady state,
the index policy is asymptotically optimal for the problem with hard capacity constraint.
\end{theorem}

We note that there is no stability in the case $\gamma=1$. This fact is in agreement with
the observed extreme unfairness of the MIMD dynamics \cite{AAP05}. However, the relaxed control in the form
of threshold policy still works in this case.

\section{Conclusion and future research}
\label{sec:conc}

We have analysed the impulsive control of the G-AIMD dynamics under hard and relaxed constraints.
In the case of the relaxed constraint, we have shown that threshold-type policy is optimal. For the
much more challenging case of hard constraint, we have proposed Whittle-type index policy and have
shown its asymptotic optimality in the homogeneous case. A number of interesting research questions
remain open: first, it would be good to establish global stability of the index policy for the
nonlinear G-AIMD. Second, we need to investigate if asymptotic optimality of the index policy
still holds in the non-homogeneous case. And third, it will be useful to extend the case of one
constraint to the case of several constraints, as in \cite{ABP17,KMT98,MM11,S12},
representing a more complex network structure.

\end{document}